\documentclass[a4paper,12pt]{article}
\usepackage{latexsym}
\usepackage{amsmath}
\usepackage{amsbsy}
\usepackage{amssymb}
\usepackage{amscd}

\begin{document}
\newcommand{\ep}{\hspace*{\fill}$\Box$}
\newcommand{\eps}{\varepsilon}
\newcommand{\pr}{{\bf Proof. }}
\newcommand{\ms}{\medskip\\}
\newcommand{\cl}{\mbox{\rm cl}}
\newcommand{\g}{\ensuremath{\mathfrak g} }
\newcommand{\gc}{${\cal G}$-complete }
\newcommand{\sa}{\stackrel{\scriptstyle s}{\approx}}
\newcommand{\prol}{\mbox{\rm pr}^{(n)}}
\newcommand{\prolo}{\mbox{\rm pr}^{(1)}}
\newcommand{\deta}{\frac{d}{d \eta}{\Big\vert}_{_{0}}}
\newcommand{\detas}{\frac{d}{d \eta}{\big\vert}_{_{0}}}
\newcommand{\R}{\mathbb R}
\newcommand{\N}{\mathbb N}
\newcommand{\C}{\mathbb C}
\newcommand{\Z}{\mathbb Z}
\newcommand{\K}{\mathbb K}
\newcommand{\sR}{\mathbb R}
\newcommand{\sN}{\mathbb N}
\newcommand{\gK}{{\cal K}}
\newcommand{\gR}{{\cal R}}
\newcommand{\gC}{{\cal C}}
\newcommand{\Dp}{${\cal D}'$ }                                          
\newcommand{\go}{${\cal G}(\Omega)$ }
\newcommand{\grn}{${\cal G}(\R^n)$ }
\newcommand{\grp}{${\cal G}(\R^p)$ }
\newcommand{\grq}{${\cal G}(\R^q)$ }
\newcommand{\gt}{${\cal G}_\tau$ }
\newcommand{\gto}{${\cal G}_\tau(\Omega)$ }
\newcommand{\gtrn}{${\cal G}_\tau(\R^n)$ }
\newcommand{\gtrp}{${\cal G}_\tau(\R^p)$ }
\newcommand{\gtrq}{${\cal G}_\tau(\R^q)$ }
\newcommand{\gtn}{{\cal G}_\tau(\R^n) }
\newtheorem{thr}{\hspace*{-1.1mm}}[section]
\newcommand{\bt}{\begin{thr} {\bf Theorem }}
\newcommand{\et}{\end{thr}}
\newcommand{\bp}{\begin{thr} {\bf Proposition }}
\newcommand{\bc}{\begin{thr} {\bf Corollary }}
\newcommand{\blem}{\begin{thr} {\bf Lemma }}
\newcommand{\bex}{\begin{thr} {\bf Example }\rm}
\newcommand{\bexs}{\begin{thr} {\bf Examples }\rm}
\newcommand{\bd}{\begin{thr} {\bf Definition }}
\newcommand{\beast}{\begin{eqnarray*}}
\newcommand{\eeast}{\end{eqnarray*}}
\newcommand{\wsc}[1]{\overline{#1}^{wsc}}
\newcommand{\todo}[1]{$\clubsuit$\ {\tt #1}\ $\clubsuit$}
\newcommand{\rem}[1]{\vadjust{\rlap{\kern\hsize\thinspace\vbox%
                       to0pt{\hbox{${}_\clubsuit${\small\tt #1}}\vss}}}}
\newcommand{\ahat}{\ensuremath{\hat{\mathcal{A}}_0(M)} }
\newcommand{\atil}{\ensuremath{\tilde{\mathcal{A}}_0(M)} }
\newcommand{\aqtil}{\ensuremath{\tilde{\mathcal{A}}_q(M)} } 
\newcommand{\amtil}{\ensuremath{\tilde{\mathcal{A}}_m(M)} } 
\newcommand{\ehat}{\ensuremath{\hat{\mathcal{E}}(M)} } 
\newcommand{\emhat}{\ensuremath{\hat{\mathcal{E}}_m(M)} }
\newcommand{\nhat}{\ensuremath{\hat{\mathcal{N}}(M)} }
\newcommand{\ghat}{\ensuremath{\hat{\mathcal{G}}(M)} } 
\newcommand{\lhat}{\ensuremath{\hat{L}_X} }                    
\newcommand{\comp}{\subset\subset}
\newcommand{\al}{\alpha}
\newcommand{\bet}{\beta} 
\newcommand{\ga}{\gamma}
\newcommand{\Om}{\Omega}\newcommand{\Ga}{\Gamma}\newcommand{\om}{\omega}
\newcommand{\si}{\sigma}\newcommand{\la}{\lambda}
\newcommand{\de}{\delta}
\newcommand{\vphi}{\varphi}\newcommand{\dl}{{\displaystyle \lim_{\eta>0}}\,}
\newcommand{\intl}{\int\limits}\newcommand{\su}{\sum\limits_{i=1}^2}
\newcommand{\D}{{\cal D}}\newcommand{\Vol}{\mbox{Vol\,}}
\newcommand{\Or}{\mbox{Or}}\newcommand{\sign}{\mbox{sign}}
\newcommand{\na}{\nabla}\newcommand{\pa}{\partial}
\newcommand{\ti}{\tilde}\newcommand{\T}{{\cal T}} \newcommand{\G}{{\cal G}}
\newcommand{\DD}{{\cal D}}\newcommand{\X}{{\cal X}}\newcommand{\E}{{\cal E}} 
\newcommand{\CC}{{\cal C}}\newcommand{\vo}{\Vol}
\newcommand{\bat}{\bar t}
\newcommand{\bx}{\bar x}
\newcommand{\by}{\bar y} \newcommand{\bz}{\bar z}\newcommand{\br}{\bar r}
\newcommand{\fr}{\frac{1}}\newcommand{\il}{\int\limits}
\newcommand{\nn}{\nonumber}
\newcommand{\supp}{\mathop{\mathrm{supp}}}

\newcommand{\vp}{\mbox{vp}\frac{1}{x}}\newcommand{\A}{{\cal A}}
\newcommand{\Ll}{L_{\mbox{\small loc}}}\newcommand{\Hl}{H_{\mbox{\small loc}}}
\newcommand{\Lll}{L_{\mbox{\scriptsize loc}}}
\newcommand{\be}{ \begin{equation} }\newcommand{\ee}{\end{equation} }
\newcommand{\beq}{ \begin{equation} }\newcommand{\eeq}{\end{equation} }
\newcommand{\bea}{\begin{eqnarray}}\newcommand{\eea}{\end{eqnarray}}
\newcommand{\beas}{\begin{eqnarray*}}\newcommand{\eeas}{\end{eqnarray*}}
\newcommand{\beqs}{\begin{equation*}}\newcommand{\eeqs}{\end{equation*}}
\newcommand{\lb}{\label}\newcommand{\rf}{\ref}
\newcommand{\GL}{\mbox{GL}}\newcommand{\bfs}{\boldsymbol}
\newcommand{\ben}{\begin{enumerate}}\newcommand{\een}{\end{enumerate}}
\newcommand{\ba}{\begin{array}}\newcommand{\ea}{\end{array}}
\newtheorem{thi}{\hspace*{-1.1mm}}[section]
\newcommand{\bthm}{\begin{thr} {\bf Theorem. }}
\newcommand{\bprop}{\begin{thr} {\bf Proposition. }}
\newcommand{\bcor}{\begin{thr} {\bf Corollary. }}
\newcommand{\bdef}{\begin{thr} {\bf Definition. }}
\newcommand{\brem}{\begin{thr} {\bf Remark. }\rm}
\newcommand{\bth}{\begin{thr}\rm}
\newcommand{\ethi}{\end{thr}}
\newcommand{\ca}{{\cal A}}
\newcommand{\cb}{{\cal B}}
\newcommand{\cc}{\mathcal{C}}
\newcommand{\cd}{{\cal D}}
\newcommand{\ce}{{\cal E}}
\newcommand{\cg}{{\cal G}}
\newcommand{\ci}{{\cal I}}
\newcommand{\cn}{{\cal N}}
\newcommand{\cs}{{\cal S}}
\newcommand{\rmd}{\mbox{\rm d}}
\newcommand{\io}{\iota}
\newcommand{\bnot}{\begin{thr} {\bf Notation }}
\newcommand{\lgl}{\langle}
\newcommand{\rgl}{\rangle}
\newcommand{\spp}{\mbox{\rm supp\,}}

\begin{center}
{\bf \Large A global theory of algebras of generalized functions}
\vskip2mm
{\large     M. Grosser, M. Kunzinger, R. Steinbauer

        Universit\"at Wien\\ 
        Institut f\"ur Mathematik\\[1mm]

          J. Vickers \\ University of Southampton, \\ Department of Mathematics}
\end{center}
\vskip3mm
{\small
{\sc Abstract.}
We present a geometric approach to defining an algebra $\hat\G(M)$ 
(the Colombeau algebra) of generalized functions on a smooth manifold $M$
containing the space $\D'(M)$ of distributions on $M$. Based 
on differential calculus in convenient vector spaces we achieve an intrinsic 
construction of $\hat\G(M)$. $\hat\G(M)$ is a 
{\em differential} algebra, its elements possessing Lie derivatives with respect 
to arbitrary smooth vector fields. Moreover, we construct a canonical linear embedding of 
$\D'(M)$ into $\hat\G(M)$ that renders $\cc^\infty (M)$ a 
faithful subalgebra of $\hat\G(M)$. Finally, it is shown that this embedding commutes 
with Lie derivatives. Thus $\hat \G(M)$ retains all the distinguishing
properties of the local theory in a global context.
}
\setcounter{page}{1}
\section{Introduction}
Colombeau's theory of algebras of generalized functions (\cite{C1}, \cite{C2},
\cite{C3}, \cite{mo}) is a well-established tool for treating
nonlinear problems involving singular objects, in particular, for
studying nonlinear differential equations in generalized functions. 
A Colombeau algebra $\mathcal{G}(\Om)$ on an open subset $\Om$ of $\R^n$ is a 
differential algebra containing $\D'(\Om)$ as a linear subspace and 
$\cc^\infty(\Om)$ as a faithful subalgebra. The embedding $\D'(\Om) \hookrightarrow
\mathcal{G}(\Om)$ commutes with partial differentiation and the functor $\Om \to
\mathcal{G}(\Om)$ defines a fine sheaf of differential algebras. In view of L.
Schwartz's impossibility result (\cite{Schw}) these properties of Colombeau algebras
are optimal in a very precise sense (cf.\  \cite{mo}, \cite{survey}).
In the so-called  ``full'' version of the theory the embedding $D'(\Om) 
\hookrightarrow \mathcal{G}(\Om)$ is canonical, as opposed to the
``special'' or ``simplified'' version (cf.\  the remark below), where
the embedding depends on a particular mollifier.

The main interest in the theory so far has come from the field 
of nonlinear partial differential equations (cf.\ e.g.\ \cite{bmo},
\cite{bmo2}, \cite{chmo}, \cite{chmo2}, \cite{cmo}, as well as \cite{mo} and
the literature cited therein), whereas the development of a theory
of Colombeau algebras on manifolds proceeded at a much slower pace. \cite{AB} 
presents an approach which basically consists in lifting the sheaf $\mathcal{G}$
from $\R^n$ to a manifold $M$. A more refined sheaf-theoretic analysis of 
Colombeau algebras on manifolds is  given  in  \cite{DDR}.
It treats the ``special'' version
of the algebra in the sense of~\cite{mo}, p.~109f, whose elements 
(termed ultrafunctions by the authors) depend on
a   real   regularization  parameter.    
In  both approaches,
as well as in the construction envisaged in~ \cite{Bal}, the
canonical embedding $\cc^\infty(M) \hookrightarrow \D'(M) \hookrightarrow
\mathcal{G}(M)$ is lost when passing from $M=\Om \subseteq \R^n$ to a general
manifold (due to the fact that in these versions action of diffeomorphisms 
commutes with embedding only in the sense of association but not with
equality in the algebra). Also, although Lie derivatives of elements of $\mathcal{G}(M)$ 
are defined in \cite{DDR} the operation of taking Lie derivatives does not commute
with the embedding (again this property only holds in the sense of association).

To remedy the first of these defects, J. F. Colombeau and A. Meril in \cite{CM}
(using earlier ideas of~\cite{C1})
introduced an algebra of generalized functions on $\Om$ whose elements depend
smoothly (in the sense of Silva-differentiability) on $(\vphi,x)$, where 
$\vphi \in \D(\R^n)$ and $x\in \Om$. The aim of \cite{CM} is
to make the embedding of $\D'(\Om)$ (which is done basically by convolution with
the  $\vphi$'s:  $u\in  \D'(\Om)  \mapsto  ((\vphi,x)\mapsto  \langle u(y),
\vphi(y-x)\rangle)$) commute with the action of diffeomorphisms.  
However, an explicit counterexample in \cite{Jel} demonstrated that the 
construction given in \cite{CM} is in fact not diffeomorphism 
invariant. Also in \cite{Jel}, J. Jel\'\i nek gave an improved version
of the theory clarifying a number of open questions but still falling
short of establishing the existence of an invariant (local) Colombeau 
algebra.
Finally, this aim was achieved in   
\cite{fo1} and \cite{fo2}, where a complete  construction of  diffeomorphism
invariant   Colombeau   algebras   on   open  subsets  of  $\R^n$ is developed.  
Several characterization  results  derived  there will be crucial for the
presentation in the following sections. 

Summing up, the current state of the theory is that a diffeomorphism invariant 
local version (and, therefore, a version on manifolds)
of Colombeau algebras providing canonical embeddings of smooth functions and
distributions is available. Recent applications of Colombeau algebras to
questions  of  general  relativity (e.g.\  \cite{CVW}, \cite{Bal}, \cite{KS},
\cite{VW}, for a survey see \cite{Vesi}) have
underscored the need for a theory of algebras of generalized functions on manifolds that 
enjoys two additional features: first, it should be {\it geometric} in the
sense that its basic objects should be defined intrinsically on the manifold
itself. Second, the object to be constructed should be a {\it
differential algebra} with 
Lie derivatives commuting with the embedding of $\D'(M)$. 

In this paper we construct an algebra $\hat\mathcal{G}(M)$ satisfying both of
these requirements. In particular, elements of $\hat\mathcal{G}(M)$ possess a
Lie derivative $\hat L_X$ with respect to arbitrary smooth vector fields $X$ 
on $M$. Moreover, each Lie derivative commutes with the canonical embedding of $\D'(M)$
into $\hat\mathcal{G}(M)$, so that in fact {\it all} the distinguishing
properties of Colombeau algebras on open subsets of $\R^n$ are retained in the
global case. The key concept leading to a global formulation of the theory is
that of {\em smoothing kernels}: Definition \ref{kernels} (i) below is in a
sense the  diffeomorphism invariant `essence' of the process of regularization 
via convolution and linear scaling on $\R^n$ while \ref{kernels} (ii) is the
invariant formulation of the interplay between $x$- and $y$-differentiation
in the local context. Finally, a number of localization results in section
\ref{localization} allow one to make full use of the well-developed local theory
also in the global context.

\section{Notation and terminology, the local theory} \label{notterm}
Throughout this paper, $M$ will denote an oriented paracompact
$\cc^\infty$-manifold of dimension $n$. An atlas of $M$ will usually be
written in the form $\mathfrak{A} = \{(U_\alpha,\psi_\alpha)\ :\
\alpha \in A\}$. By $\Omega^n_c(M)$ we mean the space of compactly
supported (smooth) $n$-forms on $M$. Locally, for coordinates $y^1,\dots,y^n$
on $U_\al$, elements of $\Om_c^n(\psi_\al(U_\al))$ will be written as
$\vphi\, d^ny := \vphi\, dy^1 \wedge \dots \wedge dy^n$.  The pullback
of any $\om\in \Omega^n_c(\psi_\al(U_\al))$ under $\psi_\al$ is
written as $\psi_\al^*(\om)$.  Then $\int_M \psi_\al^* (\vphi\,d^ny) =
\int_{\psi_\al(U_\al)} \vphi\,d^ny$ for all $\vphi\,d^ny \in
\Omega^n_c(\psi_\al(U_\al))$. For $U(\mathrm{open}) \subseteq M$ open we will
notationally suppress the embedding of $\Om^n_c(U)$ into $\Om_c^n(M)$,
and similarly for the inclusion of $\D(U)$ (compactly supported smooth
functions with support contained in $U$) into $\D(M)$.  We generally
use the following convention: if $B\subseteq \psi_\al(U_\al)$ then
$\hat B := \psi_\al^{-1}(B)$ and if $f: \psi_\al(U_\al)\to \R$ resp.\ $\C$ then
$\hat f := f\circ\psi_\al$.  A similar ``hat-convention''
will be applied to the function spaces to be defined in the following
section. Moreover, since $M$ is supposed to be oriented we can and shall
identify $n$-forms and densities henceforth.

For an open subset $\Om$ of $\R^n$, the space of distributions on
$\Om$ (i.e., the dual of the (LF)-space $\D(\Om)$) will be denoted by
$D'(\Om)$.  For a diffeomorphism $\mu: \tilde\Om \to \Om$, the
pullback of any $u\in \D'(\Om)$ under $\mu$ is defined by
\begin{equation} \label{distpull}
\langle \mu^*(u),\varphi\rangle = \langle u(y),
\varphi(\mu^{-1}(y))\cdot|\mathrm{det}\, D\mu^{-1}(y)|\rangle
\end{equation}
Concerning distributions on manifolds we follow the terminology
of~\cite{deRham} and \cite{Marsden}, so the space of distributions on
$M$ is defined as $\D'(M) = \Omega^n_c(M)'$.  Observe that in this
setting test objects no longer have {\em function} character but are
$n$-forms.  
In the context of Colombeau algebras it is natural to regard smooth
functions as regular distributions (which in the non-orientable case
enforces the use of test {\em densities}; this is in accordance 
with~\cite{Hoe} but has to be distinguished from the setting 
of~\cite{D3}). 

Operations on
distributions are defined as (sequentially) continuous extensions of
classical operations on smooth functions.
In particular, for $X\in \mathfrak{X}(M)$ (the space of smooth vector fields 
on $M$) and $u\in \D'(M)$ the
Lie derivative of $u$ with respect to $X$ is given by $\langle L_X u,
\om\rangle = -\langle u, L_X \om \rangle$. If $u\in \D'(M)$,
$(U_\al,\psi_\al) \in \mathfrak{A}$ then the local representation of
$u$ on $U_\al$ is the element $(\psi_\al^{-1})^*(u)\in
\D'(\psi_\al(U_\al))$ defined by
\begin{equation}  \label{locdist}
\langle (\psi_\al^{-1})^*(u), \vphi \rangle = \langle u|_{U_\al},
\psi_\al^*(\vphi d^n y) \rangle \qquad \forall \vphi \in
\D(\psi_\al(U_\al))
\end{equation}
It should be clear from (\ref{distpull}) and (\ref{locdist}) that the
character of test objects as $n$-forms is actually already built into the
local theory of distributions on $\R^n$ (i.e., on the right hand side 
of (\ref{distpull}) $u$ acts exactly on the coefficient 
function of $(\mu^{-1})^*(\vphi\, d^ny)$).

As in \cite{fo1}, \cite{fo2} differential calculus in infinite dimensional 
vector spaces will be
based on the presentation in \cite{KM}. The basic idea is that a map
$f:E\to F$ between locally convex spaces is smooth if it transports
smooth curves in $E$ to smooth curves in $F$, where the notion of
smooth curves is straightforward (via limits of difference
quotients). The
diffeomorphism invariant theory of Colombeau algebras on open subsets
of $\R^n$ introduced in $\cite{fo1}$ is
based on this notion of differentiability.  
Of the two (equivalent) formalisms for describing the
diffeomorphism invariant local Colombeau theory analyzed in \cite{fo1},
Section 5
only one (the so-called J-formalism which was employed by Jel\'\i nek
in \cite{Jel}) lends 
itself naturally to an intrinsic generalization on manifolds, as the
embedding of distributions does not involve a translation 
(see below). We briefly recall the main features of this theory.

Let $\Om \subseteq \R^n$ open; then define
\beas 
 \A_0(\Om)\! &=&\! \{\vphi \in \D(\Om)| \int\! \vphi(\xi)d\xi =
1\}\\ 
 \A_q(\R^n)\! &=&\! \{\vphi \in \A_0(\R^n)| \int\! \vphi(\xi)\xi^\al
d\xi = 0, \ 1\le |\al| \le q,\ \al \in \N_0^n \} \quad (q \in \N) \,.
\eeas
The basic space of the diffeomorphism invariant local Colombeau algebra 
is defined to be $\E(\Om) = \cc^\infty(\A_0(\Om)\times \Om)$.    
The algebra itself is constructed as the quotient of the space of moderate 
modulo the ideal of negligible elements $R$ of the basic space where the respective 
properties are defined by plugging  scaled and transformed ``test objects''
into $R$ and analyzing its asymptotic behavior on these ``paths'' as the 
scaling parameter $\eps$ tends to $0$.  Diffeomorphism 
invariance of the whole construction is
achieved by diffeomorphism invariance of this process--- termed ``testing
for moderateness resp.\  negligibility''---in~\cite{fo1}, Section 9. 
We shall discuss this matter in some more detail at the
end of this section but now  proceed by defining the actual ``test objects.''
Set $I=(0,1]$.
Let $\cc^\infty_b(I\times\Om,\ca_0(\R^n))$ be the space of smooth maps 
$\phi:I\times\Om\to\ca_0(\R^n)$ such that for each $K \comp \Om$
(i.e., $K$ a compact subset of $\Om$) and any $\al\in\N_0^n$, the set 
$\{\pa^\al_x\phi(\eps,x)\mid\eps\in(0,1],x\in  K\}$
is bounded in $\cd(\R^n)$. For any $m\ge 1$ we set
\beas 
&& \!\!\!\!\!\!\!\! \A_m^\Box(\Om) := \{\Phi \in
C_b^\infty(I\times\Om,\A_0(\R^n)) |\,\sup_{x\in K}|\int
\Phi(\eps,x)(\xi)\xi^\al d\xi| = O(\eps^{m}) \\ &&
\!\!\!\!\!\!\!\!\hphantom{\A_m^\Box(\Om) := \{\Phi \in
C_b^\infty(I\times\Om,\A_0(\R^n)) |} (1 \le |\al| \le m) \ \forall
K\comp \Om \} \\ 
&& \!\!\!\!\!\!\!\! \A_m^\Delta(\Om) := \{\Phi \in
C_b^\infty(I\times\Om,\A_0(\R^n)) |\,\sup_{x\in K}|\int
\Phi(\eps,x)(\xi) \xi^\al d\xi| = O(\eps^{m+1-|\al|}) \\ &&
\!\!\!\!\!\!\!\!\hphantom{\A_m^\Box(\Om) := \{\Phi \in
C_b^\infty(I\times\Om,\A_0(\R^n)) |} (1 \le |\al| \le m) \ \forall
K\comp \Om \}  \eeas Elements of $\A_m^\Box(\Om)$ are said to have
asymptotically vanishing moments of order $m$ (more precisely, in the
terminology of \cite{fo2}, elements of $\A_m^\Box(\Om)$ are of type
$[A_g]$, the abbreviation standing for {\em asymptotic} vanishing of
moments {\em globally}, i.e., on each $K\comp \Om$).  
The spaces  $\A_m^\Box(\Om)$
and $ \A_m^\Delta(\Om) $ play a crucial role in the characterizations 
of    the    algebra   (\cite{fo1},   Section 10) and in
Section~\ref{localization} below (cf.\  also the discussion of diffeomorphism
invariance at the end of this section). Also for later use we note 
that $\A_m^\Box(\Om)\subseteq \A_m^\Delta(\Om) \subseteq
\A_{2m-1}^\Box(\Om)$. For $x\in \R^n$, $\eps \in I$
we define translation resp.\   scaling operators by $T_x : \D(\R^n) \to
\D(\R^n)$, $T_x(\vphi) = \vphi(.  - x)$ and $S_\eps: \D(\R^n) \to
\D(\R^n)$, $S_\eps\vphi = \eps^{-n}\vphi(\frac{.}{\eps})$.
Now the  subspaces of moderate
resp.\  negligible elements of $\E(\Om)$ are defined by    
\beas
&& \!\!\!\! \E_m(\Om) =  \{R \in \E(\Om) \mid \forall K\comp \Om \,\forall \al \in \N_0^n
\,\exists N\in \N \, \\
&& \!\!\!\! \forall \phi \in \cc^\infty_b(I\times\Om,\ca_0(\R^n)) : 
\sup\limits_{x\in K} |\pa^\al(R(T_x S_\eps\phi(\eps,x),x))| = O(\eps^{-N}) 
\ \  (\eps\to 0)\} \\
&& \!\!\!\! \mathcal{N}(\Om) =  \{R \in \E(\Om) \mid \forall K\comp \Om \,
\forall \al \in \N_0^n\, 
\forall     r\in \N \,  \exists    m\in \N \, \\
&& \!\!\!\! \forall \phi \in \cc^\infty_b(I\times\Om,\ca_m(\R^n)) : 
\sup\limits_{x\in K} |\pa^\al(R(T_x S_\eps\phi(\eps,x),x))| = O(\eps^{r}) 
\ \ (\eps\to 0) \}
\eeas
By \cite{fo1}, Th.\ 7.9 and \cite{fo2}, Corollaries 16.8 and 17.6,
$R\in
\mathcal{E}_m(\Om)$ is negligible iff
$\forall K\subset\subset\Om\ \forall\al\in\N_0^n\ \forall r\in \N \
\exists m\in\N 
\ \forall\phi\in \ca^\Box_m(\Om),\ 
\exists C>0\ \exists \eta>0\ \forall\eps\,(0<\eps<\eta)
    \ \forall x\in K$:
\beas
|\pa^\al(R(T_x S_\eps\phi(\eps,x),x))|\leq C\eps^{r}\,.
\eeas
The  diffeomorphism  invariant Colombeau algebra $\G(\Om)$ on $\Om$ is 
the  quotient  algebra $\E_m(\Om)/\mathcal{N}(\Om)$. Partial derivatives in
$\G(\Om)$ are defined as
\begin{equation} \label{jeldiff}
(D_i R)(\vphi,x) = -((d_1 R)(\vphi,x))(\pa_i \vphi) + (\pa_i R)(\vphi,x)
\end{equation}
where  $d_1$ and $\pa_i$ denote differentiation with respect to $\vphi$ and
$x_i$, respectively. 
Note that formula~(\ref{jeldiff}) is exactly 
the result of  translating   the  usual  partial  differentiation  from  the
C-formalism  to  the  J-formalism  (cf.\ \cite{fo1},  Section 5). 
With these operations, $\Om \to \G(\Om)$ becomes a 
fine  sheaf  of  differential  algebras  on $\R^n$. For $R\in \E_m(\Om)$ we
denote by $\cl[R]$ its equivalence class in $\G(\Om)$.
The map
\beas
 \iota: \D'(\Om) & \to & \G(\Om)  \\
 u & \to & \cl[ (\vphi,x) \to \langle u , \vphi \rangle]
\eeas
provides  a linear  embedding  of  $\D'(\Om)$ into $\G(\Om)$ whose restriction to
$\cc^\infty(\Om)$ coincides with the embedding
\beas
 \sigma: \cc^\infty(\Om) & \to & \G(\Om)  \\
 f & \to & \cl[ (\vphi,x) \to f(x)]
\eeas
so $\iota$ renders $\cc^\infty(\Om)$  a faithful subalgebra and $\D'(\Om)$ a linear
subspace of $\G(\Om)$. Moreover, $\iota$ commutes with partial derivatives due to
the specific form of (\ref{jeldiff}).
The   pullback   of  $R\in  \E(\Om)$  under  a
diffeomorphism $\mu: \tilde\Om \to \Om$ is defined as
$$
(\hat\mu R)(\tilde\vphi,\tilde x) = R(\bar\mu(\tilde\vphi,\tilde x))
$$
where  $\bar\mu(\tilde\vphi,\tilde  x) = ((\tilde\vphi\circ\mu^{-1})
\cdot |\det  D  \mu^{-1}|,\mu(\tilde  x))$.  Pullback  under diffeomorphisms then
commutes with the embedding of $\D'(\Om)$ into $\G(\Om)$.

Finally we return to the issue of diffeomorphism invariance. Since the
action of a diffeomorphism on $R\in\E(\Om)$ is defined in a functorial way
diffeomorphism-invariance of the definition of the moderate resp.\  negligible 
elements of $\E(\Om)$ indeed is invariant if the class of scaled and transformed 
``test objects'' is; more precisely if 
\bea
        \phi(\eps,x)&=&
        S_\eps^{-1}\circ T_{-x}\circ\mathrm{pr}_1\circ\bar\mu(T_{\mu^{-1} x}
        \circ S_\eps\tilde\phi(\eps,\mu^{-1}x)(\xi),\mu^{-1}x)\nn\\
        &=&\tilde\phi(\eps,\mu^{-1} x)\left(\frac{\mu^{-1}(\eps\xi+x)
        -\mu^{-1}x}{\eps}\right)\cdot|\det D\mu^{-1}(\eps\xi+x)|
        \label{mubar}
\eea is a valid ``test object'' if $\tilde\phi(\eps,\tilde x)$ was.
However, if $\tilde\phi \in \cc^\infty_b(I\times\Om,\ca_q(\R^n))$ in general
$\phi$ will neither be defined on all of $I\times \Om$ nor will its moments be 
vanishing up to order $q$. To remedy these defects a rather delicate analysis
of the testing procedure is required: 
Denote by
$\cc^\infty_{b,w}(I\times\Om,\ca_0(\R^n))$ the space of all $\phi: D \to
\ca_0(\R^n)$  where $D$ is some subset (depending on $\vphi$) of
$(0,1]\times\Om$ and for $D,\vphi$ the following holds: 

For each
$L\subset\subset\Om$ there exists $\eps_0$ and a  subset $U$
of $D$ which is open in $(0,1]\times\Om$ such that
\bea 
&& (0,\eps_0]\times L\subseteq U(\subseteq D) \ \mbox{and} \ \phi \ 
\mbox{is smooth on } \ U \label{atoz1}\\
&& \{\pa^\bet\phi(\eps,x)\mid0<\eps\le\eps_0,\ x\in L\} \ \mbox{is bounded in} \
\cd(\R^n) \  \forall \,\bet\in\N_0^n \label{atoz2}
\eea
(the  subscript  $w$ signifies the  weaker  requirements  on  the  domain 
of definition  of $\phi$). 
Then by \cite{fo1}, Th.\ 10.5 $R\in \mathcal{E}(\Om)$ is moderate iff
$\forall K\subset\subset\Om\ \forall\al\in\N_0^n\ \exists N\in\N
\ \forall\phi\in \cc^\infty_{b,w}(I\times\Om,\ca_0(\R^n)),\ \phi : D\to\ca_0(\R^n)) 
\exists C>0\ \exists \eta>0\ \forall\eps\,(0<\eps<\eta)
    \ \forall x\in K$: $(\eps,x)\in D$ and
\beas
|\pa^\al(R(T_x S_\eps\phi(\eps,x),x))|\leq C\eps^{-N}
\eeas
Moreover, by \cite{fo1}, Th.\ 7.14, (\ref{mubar}) is an element of 
$\cc^\infty_{b,w}(I\times \tilde\Om,\ca_0(\R^n))$ for every $\phi \in
\cc^\infty_{b}(I\times\Om,\ca_0(\R^n))$.
The subspace of $\cc^\infty_{b,w}(I\times\Om,\ca_0(\R^n))$ consisting of those
$\phi$ whose moments up to order $m$ vanish asymptotically on each compact
subset of $\Om$ will be written as $\A_{m,w}^\Box(\Om)$ (also,
$\A_{m,w}^\Delta(\Om)$ is defined analogously). 
By the proof of \cite{fo1}, Cor.\  10.7 $R\in
\mathcal{E}_m(\Om)$ is negligible iff
$\forall K\subset\subset\Om\ \forall\al\in\N_0^n\ \forall r\in \N \
\exists m\in\N 
\ \forall\phi\in \A_{m,w}^\Box(\Om),\ \phi : D\to\ca_0(\R^n)) 
\, \exists C>0\ \exists \eta>0\ \forall\eps\,(0<\eps<\eta)
    \ \forall x\in K$: $(\eps,x)\in D$ and
\beas
|\pa^\al(R(T_x S_\eps\phi(\eps,x),x))|\leq C\eps^{r}
\eeas
and if $\phi\in \A_{m,w}^\Box(\Om)$ then (\ref{mubar}) defines an element
of $\A_{[\frac{m+1}{2}],w}^\Box(\tilde\Om)$. These facts directly imply
diffeomorphism invariance of local Colombeau algebras (\cite{fo1}, Thms. 
7.15 and 7.16).

The diffeomorphism invariance of the {\em scaled, transformed}  ``test objects''
demonstrated above suggests to choose an analogue of these as the main ``test objects'' 
on the manifold,
where no natural scaling resp.\  translation operator is available. It is 
precisely this notion which is captured in the definition of the smoothing kernels
below. 

\section{Smoothing kernels and basic function spaces}
In this section we introduce the basic definitions and operations needed for
an {\em intrinsic} definition of algebras of generalized functions on manifolds.
\bd
\[
\hat{\mathcal{A}}_0(M) := \{\om \in \Omega^n_c(M)\ : \ \int \om = 1 \} 
\]
\et
The basic space for the forthcoming definition of the Colombeau algebra on $M$ is defined
as follows:
\bd
\[
\hat{\mathcal{E}}(M) = \cc^\infty(\hat{\mathcal{A}}_0(M) \times M) 
\]
\et
Then 
\[
(\psi_\al^* \times \psi_\al^{-1})(\hat\mathcal{A}_0(\psi_\al(U_\al))\times \psi_\al(U_\al))
\subseteq \hat{\mathcal{A}}_0(M) \times M 
\]
and locally we have
$
\hat\mathcal{E}(\psi_\al(U_\al)) \cong \mathcal{E}(\psi_\al(U_\al))
$, the  
isomorphism being effected by $\hat{\mathcal{E}}(\psi_\al(U_\al))\ni R
\mapsto [(\vphi,x) \mapsto R(\vphi\ d^ny,x)]$. In what follows, we will therefore 
use $\hat{\mathcal{E}}(\psi_\al(U_\al))$ and  $\mathcal{E}(\psi_\al(U_\al))$
interchangeably. Clearly the map
\[
  \psi_\al^* \times \psi_\al^{-1}: \Omega^n_c(\psi_\al(U_\al))\times \psi_\al(U_\al)
  \to \Omega^n_c(U_\al) \times U_\al \hookrightarrow \Omega^n_c(M)\times M
\]
is smooth. Therefore, for any $R\in \ehat$ its local representation 
\[
(\psi_\al^{-1})^\wedge R := R \circ  (\psi_\al^* \times \psi_\al^{-1}) 
\]
is  an  element of $\hat\mathcal{E}(\psi_\al(U_\al))$. More generally, if $\mu:
M_1  \to  M_2$ is a diffeomorphism and $R\in \hat\E(M_2)$ then its pullback
$\hat\mu(R)  \in  \hat\E(M_1)$  under  $\mu$  is defined as $R\circ \bar \mu$
where $\bar\mu(\om,p) = ((\mu^{-1})^*\om,\mu(p))$ and clearly $(\mu_1 \circ
\mu_2)\hat{} = \hat\mu_2 \circ \hat\mu_1$ for diffeomorphisms $\mu_1$, $\mu_2$.

Let $f:M\times M \to \bigwedge^n T^*M$ be smooth such that for each $(p,q)\in M\times M$,
$f(p,q)$ belongs to the fiber over $q$ or, equivalently, that for each fixed $p\in M$,
$f_p := (q\mapsto f(p,q))$ represents a member of $\Om^n(M)$. Obviously, $p\mapsto f_p  
\in \cc^\infty(M,\Om^n(M))$ in this case. Given a smooth vector field $X$ on $M$, we
define two notions of Lie derivatives of $f$ with respect to $X$ which, essentially, 
arise as Lie derivatives of $q \mapsto f(p,q)$ resp.\ $p\mapsto f(p,q)$: 

On the one hand, viewing $\Om^n(M)$ together with the topology of pointwise (i.e., 
fiberwise) convergence as a locally convex space we define
\begin{equation}
(L'_X f)(p,q):= L_X(p \mapsto f(p,q)) = \left. \frac{d}{dt} 
\right|_0(f(\mathrm{Fl}^X_t)(p),q);
\end{equation}
on the other hand, we set
\begin{equation}
(L_X f)(p,q) := L_X(q \mapsto f(p,q)) = \left. \frac{d}{dt} 
\right|_0((\mathrm{Fl}^X_t)^*f_p)(q)
\end{equation}
where the latter symbol $L_X$ denotes the usual Lie derivative on the bundle
$\Om^n(M)$.

Now we are ready to introduce the space of smoothing kernels which will serve as an
analogue for the (unbounded) sequences $(\vphi_\eps)_\eps$ used in~\cite{C2}.
\bd \label{kernels}
$\Phi\in \cc^\infty(I \times M,\ahat)\subseteq \cc^\infty(I\times M\times M,\Lambda^nT^*M)$ 
is called a smoothing kernel if it satisfies the
following conditions
\begin{itemize}
  \item [(i)] $\forall K\subset\subset M$ $\exists\, \eps_0$, $C>0$ $\forall p\in K$
              $\forall \eps \le \eps_0$:  $\supp\Phi(\eps,p)\subseteq B_{\eps C}(p)$
  \item [(ii)] $\forall K \subset\subset M$ $\forall k,\, l\in \N_0$ 
        $\forall X_1,\dots,X_k,Y_1,\dots,Y_l\in\mathfrak{X}(M)$ 
        \[
                  \sup_{{p\in K}\atop{q\in M}} \|L_{Y_1}\dots
                  L_{Y_l}(L'_{X_1}+L_{X_1})\dots
            (L'_{X_k}+L_{X_k})
                        \Phi(\eps,p)(q) \| = O(\eps^{-(n+l)})
                \]
\end{itemize}
The space of smoothing kernels on $M$ is denoted by \atil.
\et
In (i) the radius of the ball $B_{\eps C}(p)$ has to be measured with respect to
the Riemannian distance induced by a Riemannian metric $h$ on $M$. 
By Lemma \ref{Riemann} below, (i) then holds in fact for the distance induced by 
{\em any} Riemannian metric $h'$ with a new set of constants $\eps_0(h')$ and 
$C(h')$. Similarly in (ii), $\|\,.\,\|$ denotes the norm induced on the fibers of 
$\Om_c^n(M)$ by any Riemannian metric on $M$ (i.e. convergence with respect to 
$\|\,.\,\|$ amounts to convergence of all components in every local chart). 
Thus both (i) and (ii) are independent of the Riemannian metric chosen on $M$,
hence intrinsic.

\blem \label{Riemann}
Let $M$ be a smooth paracompact manifold and let $h_1$, $h_2$ be Riemannian metrics
on M.
Then for all $K\subset\subset M$ there exist $\eps_0(K)$ and $C>0$ 
such that $\forall p\in K$ $\forall \eps \le \eps_0$:
  \[
     B_\eps^{(2)}(p) \subseteq B_{C\eps}^{(1)}(p)
  \]
  where $B^{(i)}_\eps(p) = \{q\in M \ : \ d_i(p,q) < \eps\}$ and $d_i$ denotes Riemannian 
distance with respect to $h_i$.
\et
\pr
We first show that $\forall K\subset\subset M$ $\exists C\geq 0$ such that
\begin{equation}\label{(i)}
       h_1(p)(v,v) \le C h_2(p)(v,v) \qquad \forall p\in K \, \forall v\in T_p M.
\end{equation}
Without loss of generality we may assume $K\comp U_\al$ for some chart $(\psi_\al,U_\al)$.
Denoting by $h_i^\al$ ($i=1,2$) the local representations of $h_i$ in this 
chart it follows that
\[
f: (x,v) \to \frac{h_1^\al(x)(v,v)}{h_2^\al(x)(v,v)}
\]
is continuous (even smooth) on $\psi_\al(U_\al)\times \R^n\setminus\{0\}$. Thus 
\[
\sup_{{x\in \psi_\al(K)}\atop{v\in \R^n\setminus\{0\}}} f(x,v) =
\sup_{{x\in \psi_\al(K)}\atop{v\in \pa B_1^\mathrm{eucl}(0)}} f(x,v) <\infty 
\]
with $B_1^\mathrm{eucl}(0)$ the Euclidian ball of radius $1$ around $0$.

Next we choose a geodesically convex (with respect to $h_2$) relatively
compact neighborhood $U_p$ of $p\in M$ (cf.\  e.g.\ \cite{oneill}, Prop.\  5.7). 
Moreover, let $C$ be the square root
of the constant in~(\ref{(i)}) with $K=\overline{U}_p$. Given any $q,q'\in U_p$ let
$\al$ be the unique $h_2$-geodesic in $U_p$ connecting $q$ and $q'$. Then
$d_1(q,q') \le L_1(\al)$ 
${\le} C L_2(\al) = C d_2(q,q')$ where $L_i$ denotes the length of $\al$ with
respect to $h_i$. 

Now 
for each $p\in K$ there exist $U_p$ and $C_p$ as above and we choose $\eps_p$ 
such that $B_{\eps_p}^{(2)}(p) \subseteq U_p$. Then there exist some 
$p_1,\dots,p_m$ in $K$ with $K\subseteq \bigcup_{i=1}^m
B_{\frac{\eps_{p_i}}{2}}^{(2)}(p_i) =: U$ and we set $\eps_0 :=
\min(\mathrm{dist}_2(K,\pa U),\frac{\eps_{p_1}}{2},\dots,\frac{\eps_{p_m}}{2})$
and $C:= \max_{1\le i \le m} C_{p_i}$. 

Let $p\in K$, $\eps\le \eps_0$, $q\in B_\eps^{(2)}(p)$. There exists some $i$ 
with $d_2(p,p_i)\le \frac{\eps_{p_i}}{2}$ and by construction $d_2(p,q)\le
\frac{\eps_{p_i}}{2}$, so $p,q \in B_{\eps_{p_i}}^{(2)}(p_i) \subseteq U_{p_i}$.
Hence $d_1(p,q) \le C_{p_i} d_2(p,q)$ and, finally, $q\in B^{(1)}_{C\eps}(p)$.
\ep\medskip\\

The next step is to introduce the following grading on the space of
smoothing kernels.
\bd \label{smoothingkernels}
For each $m\in\N$ we denote by $\tilde\mathcal{A}_m(M)$  
the set of all $\Phi\in \atil$ such 
that $\forall f\in \cc^\infty(M)$ and $\forall K\subset\subset M$
\[
\sup_{p\in K}  |f(p)- \int_M f(q)\Phi(\eps,p)(q) | = O(\eps^{m+1})
\]
\et
\brem
\ref{smoothingkernels}  is modelled with a view to reproducing the main
technical  ingredient  for  proving  $\iota|_{\cc^\infty} = \sigma$ (i.e. the
fact  that  the  embedding  of  distributions  into  $\G$  coincides with the
``identical'' embedding on $\cc^\infty$, cf.\  e.g., \cite{fo1}, Th.\ 7.4 (iii)) 
in the local theory.   Essentially,   the   argument  in  the  local  case  
consists  in (substitution of $y' = \frac{y-x}{\eps}$ and) Taylor
expansion  of  $\int  (f(x)  -  f(y)) \eps^{-n} \vphi(\frac{y-x}{\eps}) dy$
($\vphi\in \mathcal{A}_0$) yielding appropriate powers of $\eps$ as to establish
this term to be negligible. In fact, this argument is
at the very heart of Colombeau's construction and may be viewed as the main
technical motivation for the concrete form of the sets $\mathcal{A}_q$ and,
a fortiori, of the ideal $\mathcal{N}$ .
In  \ref{smoothingkernels}  we  turn  the  tables and {\em define} the sets
$\tilde\mathcal{A}_m(M)$  by  the  analogous  estimate. Moreover, as we
shall  see shortly (cf.\  \ref{local_kernels}), locally there is a one-to-one
correspondence between elements of $\tilde\mathcal{A}_m(M)$ and elements of
$\mathcal{A}_m^\Delta$,  so  the  two  approaches  are  in  fact equivalent
(although  only  one  of them, namely the one occurring in \ref{smoothingkernels} 
admits an intrinsic formulation on $M$). 

\et
We are now in a position to prove the nontriviality of the space of smoothing 
kernels as well as of the spaces $\tilde\mathcal{A}_m(M)$:
\blem \label{nonvoid} 
\begin{itemize}
\item[(i)] $\atil\not=\emptyset$. 
\item[(ii)] $\tilde\mathcal{A}_m(M)\not=\emptyset$ $(m\in \N)$.
\end{itemize}
\et
\pr (i)
Let $(U_\al,\psi_\al)_{\al\in A}$ be an oriented atlas of $M$ such that each $U_\al$ is  
relatively compact. Let $\{\hat\chi_\al | \al \in A\}$ be a subordinate partition of
unity and pick $\vphi \in \mathcal{A}_0(\R^n)$. 
For each $\al$, choose $\hat\chi^1_\al \in \D(U_\al)$ such that 
$\hat\chi^1_\al \equiv 1$ in an open neighborhood $W_\al \comp U_\al$ 
of $\supp\hat\chi_\al$. Let $\chi_\al:=\hat\chi_\al \circ \psi_\al^{-1}$,
$\chi^1_\al:=\hat\chi^1_\al \circ \psi_\al^{-1}$. 
Now set 
\[
\phi_\al^0(\eps,x)(y) := \eps^{-n} \vphi(\frac{y-x}{\eps}) \, .
\]
There exists some $\eps_0^\al = \eps_0^\al(\supp\chi_\al)$ in $(0,1]$ such that for each
$x\in \supp(\chi_\al)$ and each $0<\eps\le \eps_0^\al$ we have $\supp\phi_\al^0
(\eps,x) \comp \psi_\al(W_\al)$. Choose $\lambda_\al: \R \to [0,1]$ smooth, 
$\lambda_\al \equiv 1$ on $(0,\frac{\eps_0^\al}{3}]$ and $\lambda_\al \equiv 0$
on $(\frac{\eps_0^\al}{2},1]$. Finally, let $\om\in \hat\mathcal{A}_0(M)$. Then
we define our prospective smoothing kernel by
\begin{eqnarray}
&& \Phi(\eps,p)(q) := \nonumber\\ 
&& \sum_\al \hat \chi_\al(p)\left[\la_\al(\eps)\psi_\al^*\left(
\phi_\al^0(\eps,\psi_\al(p))(\,.\,) \chi^1_\al(\,.\,)d^ny\right) +
(1-\la_\al(\eps))\om\right]\!(q) \label{globkernel}
\end{eqnarray}
By construction, $\Phi$ is smooth and $\int \Phi(\eps,p) = 1$ for
all $\eps\in (0,1]$ and all $p\in M$.
For any given $K\comp M$, set $\eps_K := \min \frac{\eps_0^\al}{3}$ where
$\alpha$ ranges over the (finitely many) $\al$ with $K\cap U_\al \not= 0$;
then the terms in
(\ref{globkernel}) containing $\om$ vanish for $p\in K$ and $\eps \le \eps_K$.

In order to prove \ref{kernels} (i), let
$K\comp M$. 
Since for $p\in K$ and $\eps\le \eps_K$, $\supp \Phi(\eps,p)$ is contained 
in a finite union of supports 
of 
$$
q\mapsto [\psi_\al^*\left(
\phi_\al^0(\eps,\psi_\al(p))(\,.\,) \chi^1_\al(\,.\,)d^ny\right)](q)
$$ 
it suffices to 
consider only one term of the latter form. Clearly, there exist $C$ and 
($0<$) $\eps_0$ ($ \le \eps_K$) such that 
$\supp \phi_\al^0(\eps,x) \subseteq B_{\eps C}^{\mathrm{eucl}}(x)$ for $x\in \supp \chi_\al$
and all $\eps\le \eps_0$. Extending the pullback under $\psi_\al$ of a suitable cut-off of
the Euclidian metric on $\psi_\al(U_\al)$ to a Riemannian metric on $M$, the result
follows from \ref{Riemann}.

Concerning \ref{kernels} (ii), since $K = \bigcup_{i = 1}^m K_{\al_i}$, 
$K_{\al_i} \comp U_{\al_i}$, 
it suffices to estimate $\Phi$ on each $K_{\al_i} \times 
(0,\eps_i]$ for some $\eps_i>0$. Thus we may assume that $K\subseteq U_{\al_0}$ 
for some fixed $\al_0$.
Let $L$ be a compact neighborhood of $K$ in $U_{\al_0}$ and let 
$\eps\le \eps_L$ in what follows.
Each of the terms in (\ref{globkernel}) for which $\supp \hat\chi_\al$ does 
not intersect $K$ vanishes in some open neighborhood of $K$ and can therefore
be neglected. Each of the finitely many remaining terms 
vanishes outside $\supp \hat\chi_{\al}$. Since 
$K\subseteq (K\setminus \supp \hat\chi_\al) \cup (K\cap \overline{W}_\al)$
it is sufficient to let $p$ range over $K\cap \overline{W}_\al$
in the estimate. Now for $p\in L$, $\chi_\al^1$ can be omitted from
the corresponding term in (\ref{globkernel}).
Thus, using local coordinates of $(U_{\al_0},\psi_{\al_0})$ we have 
to estimate a finite number of Lie derivatives 
of terms of the form  
\bea \label{locLie}
&& (\mu^{-1})^*\left(\tilde y \mapsto \eps^{-n}\vphi\left(\frac{\tilde y - 
\tilde x}{\eps}\right) d^n\tilde y\right)  \nn \\
&& = y \mapsto \eps^{-n}\vphi\left(\frac{\mu^{-1}(y) - \tilde x}{\eps}\right)
\cdot \det D\mu^{-1}(y)\, d^n y\,,
\eea
each for $\tilde x\in \psi_{\al}(K\cap \overline{W}_\al)$, respectively, where
$\mu = \psi_{\al_0} \circ \psi_{\al}^{-1}$ 
and $\mu(\tilde z) = z$. Note that for $\tilde x$ in a
compact subset of $\psi_{\al}(U_{\al_0}\cap U_\al)$ and sufficiently 
small $\eps$ (say, $\eps\le \eps_1 (\le \eps_K)$) the support of 
$$
\tilde y \mapsto \eps^{-n}\vphi\left(\frac{\tilde y - 
\tilde x}{\eps}\right) d^n\tilde y
$$
is contained in $\psi_{\al}(U_{\al_0}\cap U_\al)$, rendering 
(\ref{locLie}) well-defined. Setting 
\begin{equation} \label{cmdiffeo}
\phi(\eps,x)(y) := \vphi\left(\frac{\mu^{-1}(x + \eps  y)
                                 -\mu^{-1}(x)}{\eps}\right)
                             \cdot \det D\mu^{-1}(x + \eps y)\,,
\end{equation}
the coefficient of the right hand side of (\ref{locLie}) can be written as
$T_{x}  S_\eps  \phi(\eps, x)$. By \cite{fo1}, Th.\ 7.14 and the remark
following it, for each compact set $L$ 
there exists $\eps_2 \, (\le \eps_1)$ such that  
$\{\pa_x^\beta\phi(\eps,x)(\,.\,) \mid 0<\eps\le\eps_2, x\in L\}$ is 
(defined and) bounded in $\D$ for each $\beta\in \N_0^n$.
Of course (\ref{cmdiffeo}) corresponds to (2) in \cite{CM} resp.\  to (42) 
in \cite{Jel}, the only
difference being that since we use an oriented atlas and our test objects
are forms the determinants are automatically positive. 

Since $L_X(\rho\,d^ny) = (L_X\rho) \, d^ny + \rho L_X(d^ny)$,
in order to verify \ref{kernels} (ii)
we have to estimate terms of the form
\begin{equation} \label{xx}
L_{Y_1}\dots L_{Y_{l'}}(L'_{X_1}+L_{X_1})\dots (L'_{X_{k'}}+L_{X_{k'}})
T_xS_\eps\Phi(\eps,x)(y)
\end{equation}
where
\[
X_i = \sum_{r_i=1}^n a^i_{r_i} \pa_{r_i} \qquad Y_j = 
\sum_{s_j=1}^n b^j_{s_j} \pa_{s_j}
\]
are local representations in $\psi_{\al_i}(U_{\al_i})$ of vector fields on $M$
and $0\le k'\le k$, $0\le l'\le l$.  
Explicitly, (\ref{xx}) is given by
\[
\prod_{j=1}^{l'}\left(\sum_{s_j=1}^n b^j_{s_j}(y)\frac{\pa}{\pa y^{s_j}}\right)
\prod_{i=1}^{k'}\left(\sum_{r_i=1}^n \left(a^i_{r_i}(x)\frac{\pa}{\pa x^{r_i}} + 
a^i_{r_i}(y)\frac{\pa}{\pa y^{r_i}}\right)\right) T_xS_\eps\phi(\eps,x)(y).
\]
Note that $\frac{\pa}{\pa y^{s_j}} T_x S_\eps \phi = \eps^{-1} T_x S_\eps 
\frac{\pa}{\pa y^{s_j}} \phi$ and
\beas
&& \left(a^i_{r_i}(x)\frac{\pa}{\pa x^{r_i}} + 
a^i_{r_i}(y)\frac{\pa}{\pa y^{r_i}}\right)T_x S_\eps \phi = \\
&& T_x S_\eps \left(a^i_{r_i}(x)\frac{\pa}{\pa x^{r_i}}
+ \frac{1}{\eps}(a^i_{r_i}(x+\eps y) - a^i_{r_i}(x))
\frac{\pa}{\pa y^{r_i}}\right) \phi\,.
\eeas
Each of the maps
\begin{equation} \label{diffqu}
(\eps,x,y) \to
  \begin{cases}\displaystyle
    \frac{a^i_{r_i}(x+\eps y) - a^i_{r_i}(x)}{\eps} &  \text{for} \quad \eps\not=0, \\
    Da^i_{r_i}(x)y & \text{for} \quad \eps= 0
  \end{cases}
\end{equation}
is  smooth  (hence  uniformly  bounded)  on  each relatively compact 
subset of its domain of definition. Due to the boundedness of 
$\{\pa_x^\bet \phi(\eps,x)(.)\mid 0<\eps\le \eps_2, x\in\psi_{\al_0}(K\cap 
\overline{W}_\al)\}$ only values of $y$ from a bounded region of $\R^n$
are relevant in (\ref{diffqu}). Thus one further restriction of the range
of $\eps$ establishes the claim.

(ii) Choose the $\phi_\al^0$ as in (i), yet this time additionally 
requiring $\vphi \in \A_m(\R^n)$. To estimate $\sup_{p\in K}|\int_M
\Phi(\eps,p)(q) \hat f(q) - \hat f(p)|$ for some $K\comp M$ and $\hat f
\in \cc^\infty(M)$ again only finitely many terms of (\ref{globkernel})
have to be taken into account; also, for small $\eps$, the terms involving
$\om$ vanish and $\chi_\al^1$ can be neglected. Then
\begin{eqnarray*}
  && \sup_{p\in K\cap \supp \hat\chi_\al} |\int_M
\psi_\al^*(\phi_\al^0(\eps,\psi_\al(p))(.)\,d^ny)(q)\hat f(q)  - \hat f(p)| =\\
&& = \sup_{x\in \psi_\al(K\cap \supp \hat\chi_\al)} |\int_{\psi_\al(U_\al)}
\frac{1}{\eps^n}\vphi(\frac{y-x}{\eps})f(y)\,d^ny  -  f(x)| = O(\eps^{m+1})
\end{eqnarray*}

\ep

Next, we introduce the appropriate notion of Lie derivative for elements of 
$\hat{\mathcal E}(M)$.
\bd For any $R\in \ehat$ and any $X\in \mathfrak{X}(M)$ we set
\begin{equation} \label{liealg}
(\hat{L}_X R)(\om,p) := - d_1R(\om,p)(L_X \om) + L_X(R(\om,\,.\,))|_p
\end{equation}
\et
Here, $d_1 R(\om,x)$ denotes the derivative of $\om \to R(\om,x)$ in the sense
of \cite{fo1}, section 4. In order to obtain a structural description of this
Lie derivative (which, at the same time, entails $\hat L_X R \in \hat \ce(M)$
for $R\in \hat \ce(M)$), for any $X \in \mathfrak{X}(M)$ define $X_{\A} \in 
\mathfrak{X}(\hat\A_0(M))$ by
\beas
 X_{\A}: \hat\A_0(M) &\to& \hat \A_{00}(M)\\
 X_{\A}(\om) &=& - L_X\om
\eeas
(where $\hat\A_{00}(M)$ is the linear subspace of $\Om_c^n(M)$ 
parallel to $\hat\A_0(M)$,
i.e. $\hat\A_{00}(M)$ $=$ $\{\om\in \Om_c^n(M)|\int \om = 0\}$). Then $X_{\A}R(\,.\,,p)(\om)
= d_1R(\om,p)(-L_X\om)$, so that
\[
(\hat L_X R)(\om,p) = (L_{(X_{\A},X)}R)(\om,p) = (L_{(-L_X\,.\,,X)}R)(\om,p)
\]
i.e. $\hat L_X$ is the Lie derivative of $R$ with respect to the smooth vector 
field $(X_{\A},X)$
on $\hat \A_0(M)\times M$. Thus, indeed $\hat L_X R \in \ehat$.
\brem \label{liederprop}
({\em Local description of $\hat L_X$}) 
Let $X\in \mathfrak{X}(M)$ and set $X_\alpha = (\psi_\al^{-1})^*(X|_{U_\al})$.
Since the derivative of the (linear and continuous) map $\psi_\al^*: 
\Om_c^n(\psi_\al(U_\al)) \to
\Om_c^n(U_\al)$ in any point equals 
the map itself, on $U_\al$ we obtain (writing $R$ in place of $R|_{U_\al}$): 
\begin{eqnarray*}
&& ((\psi_\al^{-1})^\wedge(\hat L_X R))(\vphi\,d^ny,x) = 
(\hat L_X R)(\psi_\al^*(\vphi\,d^ny),\psi_\al^{-1}(x))\\
&& = L_X(R(\psi_\al^*(\vphi\,d^ny),\,.\,))(\psi_\al^{-1}(x)) -\\ 
&& (d_1R)(\psi_\al^*(\vphi\,d^ny),\psi_\al^{-1}(x))
(\underbrace{L_X(\psi_\al^*(\vphi\,d^ny))}_{\psi_\al^*(L_{X_\alpha}(\vphi\,d^ny))}) \\
&& = [L_{X_\alpha}((\psi_\al^{-1})^\wedge R))](\vphi\,d^ny,x) - 
[d_1((\psi_\al^{-1})^\wedge R)(\vphi\,d^ny,x)](L_{X_\alpha}(\vphi\,d^ny))
\end{eqnarray*}
In particular, if $X_\al=\pa_{y^i}$ ($1\le i \le n$) then this exactly reproduces the local
algebra derivative with respect to $y^i$ given in (\ref{jeldiff}).
\et
Finally we are in a position to define the subspaces of moderate and negligible
elements of $\ehat$.  
\bd
$R\in \ehat$ is moderate if
 $\forall K\subset\subset M$  $\forall k\in \N_0$  $\exists N\in \N$
  $\forall\ X_1,\dots,X_k\in
                \mathfrak{X}(M)$ 
  $\forall\ \Phi \in \atil$ 
\begin{equation} \label{mod}
\sup_{p\in K} |L_{X_1}\dots L_{X_k}(R(\Phi(\eps,p),p)) | = O(\eps^{-N})
\end{equation}
The subset of moderate elements of \ehat is denoted by \emhat.
\et

\bd
$R\in \ehat$ is called negligible if it satisfies
\\
\parbox{12.7cm}{
\beas 
    && \forall K\subset\subset M \ \forall k, l \in \N_0 \ \exists m \in \N
    \ \forall \  X_1,\dots,X_k\in
                \mathfrak{X}(M)  \ \forall \Phi \in  \tilde\mathcal{A}_m(M) \\
    &&  \sup_{p\in K} |L_{X_1}\dots L_{X_k}(R(\Phi(\eps,p),p)) | = O(\eps^{l})          
\eeas
} \hfill \parbox{8mm}{\bea \label{*}\eea}\\ 
The set of negligible elements of \ehat will be denoted by \nhat.
\et

\section{Construction of the algebra, localization} \label{localization}
The following result is immediate from the definitions:
\bt
\begin{itemize}
\item[(i)] \emhat is a subalgebra of \ehat. 
\item[(ii)] \nhat is an ideal in \emhat.
\end{itemize}
\ep
\et

For the further development of the theory it is essential to achieve 
descriptions of both moderateness and negligibility that relate these concepts
to their local analogues~\cite{Jel,fo1}, thereby making available the host of local results 
already established for open subsets of $\R^n$ also
in the global context. As a first step, we are going to examine localization
properties of smoothing kernels:

\blem\label{local_kernels} 
Denote by $(U_\al,\psi_\al)$ a chart in $M$. 

{\em (A) Transforming smoothing kernels to local test objects.}
\begin{itemize}
\item[(i)] Let $\Phi$ be a smoothing kernel. Then the map $\phi$ 
defined by
\begin{equation} \label{lockeri}
\phi(\eps,x)(y)d^ny:=\eps^n\,((\psi_\al^{-1})^*\Phi(\eps,\psi_\al^{-1}x))
    (\eps y+x) \quad (x\in\psi_\al(U_\al),\,y\in\R^n)
\end{equation}
is an element of $\cc^\infty_{b,w}(I\times\psi_\al(U_\al),\ca_0(\R^n))$.

\item[(ii)] If, in addition, $\Phi \in \tilde{\mathcal{A}}_{m}(M)$ for some 
$m\in \N$ then $\phi\in{\mathcal{A}}^{\Delta}_{m,w}(\psi_\al(U_\al))$, i.e., 
\begin{equation} \label{n2}
\int \phi(\eps,x)(y)y^\bet dy = \mathcal{O}(\eps^{m+1-|\bet|})
\quad (1\le |\bet| \le m)
\end{equation}
uniformly on compact sets. In particular, if $\Phi \in 
\tilde{\mathcal{A}}_{2m-1}(M)$  then $\phi \in
\mathcal{A}_{m,w}^{\Box}(\psi_\al(U_\al))$.
\end{itemize}

{\em (B) Transporting local test objects onto the manifold.}
\begin{itemize}
\item[(i)] Let $\phi \in \cc_b^\infty(I\times\Om,\ca_0(\R^n))$ 
and $\Phi_1\in\atil$. Let $K\subset\subset\psi_\al(U_\al)$, 
$\chi,\chi_1\in\D(\psi_\al(U_\al))$ with $\chi\equiv 1$ on an open neighborhood 
of $K$ and $\chi_1\equiv 1$ on an open neighborhood of
$\supp\chi$.
\bea
\Phi(\eps,p)&:=&
    (1-\hat\chi(p)\la(\eps))\Phi_1(\eps,p)\nn \\
    &&+\hat\chi(p)\la(\eps)\psi_\al^*\left(\frac{1}{\eps^n}\phi(\eps,\psi_\al p)
    (\frac{y-\psi_\al p}{\eps})\chi_1(y)d^ny\right) \label{lockerii}
\eea
is a smoothing kernel (the smooth cut-off function $\la$ is defined in the proof). 
\item[(ii)] If, in addition,
$\phi\in{\mathcal{A}}^{\Delta}_{m}(\psi_\al(U_\al))$
and $\Phi_1\in \amtil$ then $\Phi\in \amtil$. 
In particular, if 
$\phi \in \mathcal{A}_m^{\Box}(\psi_\al(U_\al))$ and $\Phi_1\in \amtil$ then
$\Phi\in \amtil$.
\end{itemize}
\et
\pr 
\begin{itemize}
\item[(A)]
(i) Let $\hat D_1 := \{(\eps,p)\in I\times U_\al \, :\, \supp \Phi(\eps,p) 
\subseteq U_\al\}$ and set $D = \mathrm{int}((\mathrm{id} \times \psi_\al)(\hat D_1))$.
Then evidently $\phi$ is smooth on $D$. Furthermore (\ref{atoz1}) is obvious from (i) 
in~\ref{kernels}. Concerning (\ref{atoz2}), set
\[
\phi_0(\eps,x)d^ny := (\psi_\al^{-1})^*(\Phi(\eps,\psi_\al^{-1}(x))) \qquad
((\eps,x) \in D_1)
\]
Then we have
\begin{equation}\label{tauschen}
    \pa_y^\bet\pa_x^\ga\phi=\pa_y^\bet\pa_x^\ga S_\eps^{-1}T_x^{-1}\phi_0=
    \eps^{|\bet|}S_\eps^{-1}T_x^{-1}\pa_y^\bet(\pa_x+\pa_y)^\ga\phi_0.
\end{equation}
\noindent Let now $K\subset\subset\psi_\al(U_\al)$, $x\in K$ and $\eps\leq\eps(K)$. 
By Definition~\ref{kernels} (i) and Lemma \ref{Riemann} there
exists a constant $C$ such that $\mathrm{diam}(\supp\phi_0(\eps,x))\leq\eps C$
for all $x\in K$. But then $\mathrm{diam}(\supp\phi(\eps,x))\leq C$ for
all such $x$. Hence $\bigcup_{\gamma}\bigcup_{x\in K}
\supp\pa_x^\gamma\phi(\eps,x)$ is bounded uniformly in $\eps$. 
Moreover, observing that $L_{\pa_y}^\bet (L_{\pa_x} + L_{\pa_y})^\gamma
(\phi_0(\eps,x)(y) d^n y) = (L_{\pa_y}^\bet (L_{\pa_x} + L_{\pa_y})^\gamma
\phi_0(\eps,x)(y)) d^n y$, from Definition~\ref{kernels} (ii) we obtain 
$\pa_y^\bet(\pa_x+\pa_y)^\ga\phi_0=\mathrm{O}(\eps^{-(|\bet|+n)})$ which together with 
equation~(\ref{tauschen}) gives the desired boundedness property.

(ii) Now suppose that $\Phi \in \amtil$ and for each $|\bet|\le m$ choose
$f_\bet\in \D(U_\al)\subseteq \D(M)$ such that $f_\bet\circ\psi_\al^{-1}(x) = 
x^\bet$ in a neighborhood of $K\comp \psi_\al(U_\al)$. Then by assumption we have
\[
\sup_{x\in K}|\int_M f_\bet(q) \Phi(\eps,\psi_\al^{-1}(x))(q) - 
f_\bet(\psi_\al^{-1}(x))| = \mathcal{O}(\eps^{m+1})
\]
For $\eps$ sufficiently small this implies
\[
\sup_{x\in K}|\int_{\R^n} \eps^{-n}\phi(\eps,x)(\frac{y-x}{\eps}) y^\bet dy - 
x^\bet| = \mathcal{O}(\eps^{m+1})
\]
or, upon substituting $z = \eps^{-1}(y-x)$ and expanding:
\[
\sup_{x\in K}\left|\sum_{0<\gamma\le\bet} \left(\bet\atop\gamma\right)
\eps^{|\ga|} x^{\bet-\ga}\int_{\R^n} \phi(\eps,x)(z) z^\ga dz\right| 
= \mathcal{O}(\eps^{m+1})
\] From this, we prove (\ref{n2}) by induction with respect to $|\bet|$: 
for $\bet = e_k$ (the $k$-th unit vector) 
we obtain $|\int_{\R^n} \phi(\eps,x)(z)z_k dz|$ $=$ $\mathcal{O}(\eps^m)$.
Supposing now that (\ref{n2}) has already been proved for $|\bet|-1$,
the result for $|\bet|$ follows from
\beas
&& \sup_{x\in K}|\sum_{0<\gamma<\bet} \left(\bet\atop\gamma\right)
\eps^{|\ga|} x^{\bet-\ga}\underbrace{\int_{\R^n} \phi(\eps,x)(z)z^\ga 
dz}_{\mathcal{O}(\eps^{m+1-|\ga|})}  + \eps^{|\bet|}\int_{\R^n}
\phi(\eps,x)(z)z^\bet dz| \\
&& = \mathcal{O}(\eps^{m+1})
\eeas

\item[(B)]
(i) First note that by the boundedness assumption all supports of 
$\phi(I\times\supp(\chi))$ are in a fixed compact subset of $\R^n$. 
By \cite{fo1}, Lemma 6.2 we conclude that there exists $\eta>0$ such that for all
$\eps\leq\eta$ and $x\in\supp\chi$ $\supp\left(\frac{1}{\eps^n}\phi(\eps,x)
(\frac{.-x}{\eps})\right)\subseteq\{\chi_1\equiv1\}$. 
Now we are in the position to define $\la$ as follows.
Let $\la\in\cc^\infty(\R)$, $0\le\la\le1$ and $\la\equiv 1$ on
$(-\infty,\eta/3)$ and $\la\equiv 0$ on $(\eta/2,\infty)$. This actually
implies $\int\Phi(\eps,p)=1$ $\forall
p,\eps$. Smoothness again is evident while condition~\ref{kernels} (i) is
an easy consequence of \cite{fo1}, Lemma 6.2. \ref{kernels} (ii) then
follows exactly as (\ref{xx}) is proved in Lemma \ref{nonvoid}.

(ii) Let $\hat f\in \cc^\infty(M)$ and $p\in \hat K \comp U_\al$;
then
\begin{eqnarray*}
&&\!\!\!\! |\int_M \Phi(\eps,p)(q)\hat f(q) - \hat f(p)| \le \\
&&\!\!\!\! |1 - \hat\chi(p)\la(\eps)|\,
|\!\int_M \Phi_1(\eps,p)(q) \hat f(q)- \hat f(p)| + \\
&&\!\!\!\! |\hat\chi(p)\la(\eps)|\,|\!\int_M 
\psi_\al^*(\eps^{-n}\phi(\eps,\psi_\al(p))
(\frac{y-\psi_\al(p)}{\eps})\chi_1(y)d^ny)(q)\hat f(q)\! -\! \hat f(p)|
\end{eqnarray*}
The first summand is of order $\eps^{m+1}$ by assumption and the second one 
(apart from $\hat\chi$ and $\lambda$) for sufficiently small $\eps$ and 
setting $f=\hat f\circ\psi_\al^{-1}$ equals 
\[
\int_{\R^n} \phi(\eps,x)(y) f(x+\eps y) d^n y - f(x) = 
\!\!\!\! \sum_{0<|\bet|\le m} \frac{\pa^\bet f(x)}{\bet!} \eps^{|\bet|}
\underbrace{\int \phi(\eps,x)(y) y^\bet dy}_{\mathrm{O}(\eps^{m+1-|\bet|})}
\]
so the claim follows.
\end{itemize}
\ep

Restriction of an element $R$ of $\mathcal{E}(M)$ to an open subset $U$ of 
$M$ is defined by $R|_U := R|_{\hat\ca_0(U)\times U}$.
\bt \label{modloc}
(Localization of moderateness) Let $R\in \ehat$. Then 
\[
R\in \emhat \Leftrightarrow (\psi_\al^{-1})^\wedge (R|_{U_\al}) \in
\mathcal{E}_m(\psi_\al(U_\al)) \ \forall \al.
\]
\et
\pr
($\Rightarrow$) For $R\in \emhat$ and for $(U_\al,\psi_\al)$ 
some chart in $M$ let
$R' := (\psi_\al^{-1})^\wedge (R|_{U_\al})$. 
Let $K\comp \psi_\al(U_\al)=:V_\al$, $\bet \in \N_0^n$ 
and let $\phi\in$ $\cc^\infty_b(I\times\psi_\al(U_\al),$ $\ca_0(\R^n))$. 
We have to show that
\begin{equation}\label{modloc1}
\sup_{x\in K}|\pa^\bet (R'(T_xS_\eps\phi(\eps,x),x))| = 
\mathcal{O}(\eps^{-N})
\end{equation} 
for some $N\in \N$. To this end we fix some 
$\Phi_1\in\atil$ and define $\Phi\in \atil$ by (\ref{lockerii}). 
Let $\pa^\beta = \pa_{i_1}\dots \pa_{i_{|\bet|}}$ and for
$1\le i_j \le n$ choose $X_{i_j}\in \mathfrak{X}(M)$ such that the local expression
of $X_{i_j}$ coincides with $\pa_{i_j}$ on a neighborhood of $K$. Then for $\eps$ 
sufficiently small (\ref{modloc1}) equals
\[
\sup_{p\in  \hat K}|L_{X_{i_1}}\dots L_{X_{i_{|\bet|}}} R(\Phi(\eps,p),p)| \, ,
\]
so we are done.

($\Leftarrow$) Let  $X_1,\dots,X_k\in \mathfrak{X}(M)$ and (without loss of generality)
$\hat K \comp U_\al$ for some chart $(U_\al,\psi_\al)$. Let $\Phi\in \atil$ and
define $\phi$ by (\ref{lockeri}). Since $\phi: D \to \ca_0(\R^n)$ belongs to
$\cc^\infty_{b,w}(I\times\psi_\al(U_\al),\ca_0(\R^n))$ it follows from
\cite{fo1}, Th.\ 10.5 that given $\beta\in \N_0^n$ there exists some $N\in
\N$ and some $\eps_0>0$ such that for $x\in K$ and $\eps\le \eps_0$ we have
$(\eps,x)\in D$ and $|\pa^\beta (\psi_\al^{-1})^\wedge
(R|_{U_\al})(T_xS_\eps\phi(\eps,x),x))| = 
\mathcal{O}(\eps^{-N})$. Inserting this into the local representation of (\ref{mod})
immediately gives the result. \ep

\bt \label{locneg}
(Localization of negligibility) Let $R \in \emhat$. Then
\[
R\in \nhat \Leftrightarrow (\psi_\al^{-1})^\wedge (R|_{U_\al}) \in
\mathcal{N}(\psi_\al(U_\al)) 
\quad \forall \al.
\]
\et
\pr
($\Rightarrow$) Let $R\in \emhat$, $(U_\al,\psi_\al)$ some chart in $M$ and set
$R' := (\psi_\al^{-1})^\wedge (R|_{U_\al})$. Let $K\comp \psi_\al(U_\al)=:V_\al$, $l\in \N$ and 
$\bet \in \N_0^n$.
Set $k=|\bet|$ and choose $m\in \N_0$  such that 
\begin{equation} \label{locneg1} 
   \sup_{p\in \hat K} |L_{X_1}\dots L_{X_k}(R(\Phi(\eps,x),x)) | = O(\eps^{l})          
\end{equation}
for all $X_1,\dots,X_k \in \mathfrak{X}(M)$ and all $\Phi \in
\tilde\mathcal{A}_m(M)$. Now let $\phi\in \mathcal{A}_m^{\Box}(V_\al)$ and
construct $\Phi$ from $\phi$ according to (\ref{lockerii}). By
\ref{local_kernels}, $\Phi\in \tilde\mathcal{A}_m(M)$, so (with
$X_{i_j}$ as in the proof of \ref{modloc})
\[
\sup_{x\in K}|\pa^\bet (R'(T_xS_\eps\phi(\eps,x),x))| =
\sup_{p\in \hat K} |L_{X_{i_1}}\dots L_{X_{i_{|\bet|}}}(R(\Phi(\eps,x),x))|
= \mathcal{O}(\eps^l)
\] 
Thus the claim follows from the characterization of negligibility following the
definition of $\cn(\Om)$ (section \ref{notterm}).

($\Leftarrow$) Let $k\in \N_0$, $l\in \N$, $X_1,\dots,X_k\in \mathfrak{X}(M)$
and $\hat K \comp U_\al$. By the discussion at the end of section \ref{notterm} 
there exists $m'$ such that (with $R' = (\psi_\al^{-1})^\wedge R|_{U_\al}$)
\begin{equation} \label{locneg2}
\sup_{x\in K}|\pa^\bet (R'(T_xS_\eps\phi(\eps,x),x))| = 
\mathcal{O}(\eps^l)
\end{equation}
for all $|\bet|\le k$ and all $\phi \in \mathcal{A}_{m',w}^{\Box}(V_\al)$. Now set
$m = 2m'-1$ and let $\Phi \in \tilde \mathcal{A}_m(M)$. Then $\phi$ defined
by (\ref{lockeri}) is in $\mathcal{A}_{m',w}^{\Box}(V_\al)$ by
\ref{local_kernels} (A) (ii).
Hence inserting local representations of $X_1,\dots,X_k$, (\ref{locneg2}) 
immediately implies the validity of (\ref{locneg1}), thereby finishing the
proof. \ep\medskip\\
It was shown in \cite{fo2}, sec.\ 13 that for all variants of (local)
Colombeau 
algebras membership of any element $R$ of $\mathcal{E}_m$ to the ideal $\mathcal{N}$
can be tested on the function $R$ itself, without taking into account any
derivatives of $R$. As a first important consequence of the above localization 
results we note that this rather surprising simplification also holds true
for the global theory:
\bc Let $R\in \hat\mathcal{E}_m(M)$. Then
\[
R\in \nhat \Leftrightarrow (\ref{*}) \ \mathrm{ holds \ for } \ k=0.
\]
\et
\pr This follows directly from Th.\ \ref{locneg} by taking into account
\cite{fo1}, Th.\ 7.13 and \cite{fo2}, Th.\ 13.1. \ep\medskip\\
Moreover, stability of \emhat and \nhat under Lie derivatives also follows
from the local description: 
\bt Let $X\in \mathfrak{X}(M)$. Then
\begin{itemize}
\item[(i)] $\lhat \emhat \subseteq \emhat$.
\item[(ii)] $\lhat \nhat \subseteq \nhat$.
\end{itemize}
\et
\pr
Let $R\in \emhat$, $X\in \mathfrak{X}(M)$. By Th.\ \ref{modloc} for any chart
$(U_\al,\psi_\al)$ we have $(\psi_\al^{-1})^\wedge 
(R|_{U_\al}) \in \mathcal{E}_m(\psi_\al(U_\al))$.
Thus by \cite{fo1}, Th.\ 7.10 also $L_{X_\al}(\psi_\al^{-1})^\wedge$ 
$(R|_{U_\al})=
(\psi_\al^{-1})^\wedge (\lhat R|_{U_\al}) \in \mathcal{E}_m(\psi_\al(U_\al))$
(where $X_\al$ denotes the local representation of $X$), which, 
again by Th.\ \ref{modloc} gives the result. 
The claim for \nhat follows analogously from \cite{fo1}, Th.\ 7.11.
\ep\medskip\\
Finally we are in a position to define our main object of interest:
\bd 
\[
\ghat := \emhat / \nhat
\]
is called the Colombeau algebra on $M$.
\et
By construction, every \lhat induces a Lie derivative (again denoted by \lhat) on \ghat, so
\ghat becomes a differential algebra. If $R \in \ehat$, its class in $\hat\G(M)$ will
be denoted by $\cl[R]$.
\bt
\ghat is a fine sheaf of differential algebras on $M$.
\et
\pr This is a straightforward consequence of \cite{fo1}, Th.\ 8.1.\ep

\section{Embedding of distributions and smooth \\ functions} 
In this section we show that in the global context $\ghat$ displays 
the same set of (optimal) embedding properties as the local versions 
do on open sets of $\R^n$. Most importantly, we shall see that taking Lie
derivatives with respect to arbitrary smooth vector fields commutes with the
embedding.

To begin with,
let $u\in \D'(M)$. The natural candidate for the
image of $u$ in \ghat is $R_u(\om,x) = \langle u,\om\rangle$. We first show that
$R_u \in \emhat$ using Th.\ \ref{modloc}. Let $\om \in \D(\psi_\al(U_\al))$; then 
\begin{eqnarray*}
&&((\psi_\al^{-1})^\wedge (R_u|_{U_\al}))(\om,x) = (R_u|_{U_\al})(\psi_\al^*(\om\,
d^ny),\psi_\al^{-1}(x)) = \\
&& \langle u, \psi_\al^*(\om\, d^ny)\rangle = 
\langle (\psi_\al^{-1})^*(u|_{U_\al}),\om\rangle
\end{eqnarray*} 
Since $(\psi_\al^{-1})^*(u|_{U_\al}) \in \D'(\psi_\al(U_\al))$ it follows from
the local theory that 
indeed $(\psi_\al^{-1})^\wedge (R_u|_{U_\al}) \in \mathcal{E}_m(\psi_\al(U_\al))$.
Suppose now that $R_u \in \nhat$. By the same reasoning as above  
$(\om,x) \to \langle(\psi_\al^{-1})^*(u|_{U_\al}),\om\rangle \in \mathcal{N}(\psi_\al(U_\al))$ 
for each $\al$. Thus again by the respective local result
$(\psi_\al^{-1})^*(u|_{U_\al})=0$ for each $\al$, i.e. $u=0$. Therefore
\begin{eqnarray*}
&& \iota: \D'(M) \to \ghat\\
&& \iota(u) = \cl[(\om,x) \to \langle u,\om\rangle]
\end{eqnarray*}
is a linear embedding. What is more, as a direct consequence of (\ref{liealg})
(noting that distributions are linear and continuous, hence equal to their
differential in any point) we obtain   
\beas
&& \iota(L_X u)(\om,p) = \iota((\om,p) \to -\langle u,L_X\om\rangle) =
-d_1 R_u(\om,p)(L_X\om)  \\ 
&& + \underbrace{L_X(R_u(\om,\,.\,))|_p}_{=0} 
 = (\hat L_X R_u)(\om,p) = \hat L_X(\iota(u))(\om,p)
\eeas
i.e., $\iota$ commutes with arbitrary Lie derivatives.

The natural operation for embedding smooth functions into \ghat is given by
\begin{eqnarray*}
&& \sigma: \cc^\infty(M) \to \ghat\\
&& \sigma(f) = \cl[(\om,x) \to f(x)]
\end{eqnarray*}
Obviously, $\sigma$ is an injective algebra homomorphism that commutes with
Lie derivatives by (\ref{liealg}). Moreover, $\iota$ coincides with
$\sigma$ on $\cc^\infty(M)$. 
Making use of Th.\ \ref{locneg} this again follows directly from the local result. 
Summing up, we have 
\bt
$\iota: \D'(M) \to \ghat$, is a linear embedding that commutes
with Lie derivatives and coincides with $\sigma: \cc^\infty(M)\to \ghat$ 
on $\cc^\infty(M)$. Thus $\iota$ renders $\D'(M)$ a linear
subspace  and $\cc^\infty(M)$ a faithful subalgebra of \ghat.  
\et

The following commutative diagram illustrates the compatibility
properties of Lie derivatives with respect to embeddings
established in this section:

\vskip18pt

\begin{center}
\setlength{\unitlength}{1500sp}%
\begingroup\makeatletter\ifx\SetFigFont\undefined%
\gdef\SetFigFont#1#2#3#4#5{%
  \reset@font\fontsize{#1}{#2pt}%
  \fontfamily{#3}\fontseries{#4}\fontshape{#5}%
  \selectfont}%
\fi\endgroup%
\begin{picture}(8550,6747)(1201,-6973)
\thinlines
\put(2401,-961){\vector( 1, 0){6000}}
\put(1801,-1561){\vector( 0,-1){4800}}
\put(9001,-1486){\vector( 0,-1){4875}}
\put(2401,-6961){\vector( 1, 0){6000}}
\put(2401,-1561){\vector( 2,-3){1200}}
\put(8401,-1561){\vector(-2,-3){1165.385}}
\put(4201,-3961){\vector( 1, 0){2400}}
\put(3601,-4561){\vector(-2,-3){1200}}
\put(7201,-4561){\vector( 2,-3){1200}}
\put(1501,-961){\makebox(0,0)[lb]{\smash{\SetFigFont{12}{14.4}{\rmdefault}
{\mddefault}{\updefault}$\cc^\infty$}}}
\put(8901,-961){\makebox(0,0)[lb]{\smash{\SetFigFont{12}{14.4}{\rmdefault}
{\mddefault}{\updefault}$\cc^\infty$}}}
\put(5001,-711){\makebox(0,0)[lb]{\smash{\SetFigFont{12}{14.4}{\rmdefault}
{\mddefault}{\updefault}$L_X$}}}
\put(5001,-3661){\makebox(0,0)[lb]{\smash{\SetFigFont{12}{14.4}{\rmdefault}
{\mddefault}{\updefault}$L_X$}}}
\put(5001,-6661){\makebox(0,0)[lb]{\smash{\SetFigFont{12}{14.4}{\rmdefault}
{\mddefault}{\updefault}$\hat L_X$}}}
\put(1651,-6961){\makebox(0,0)[lb]{\smash{\SetFigFont{12}{14.4}{\rmdefault}
{\mddefault}{\updefault}$\hat\G$}}}
\put(8901,-6961){\makebox(0,0)[lb]{\smash{\SetFigFont{12}{14.4}{\rmdefault}
{\mddefault}{\updefault}$\hat\G$}}}
\put(3363,-5561){\makebox(0,0)[lb]{\smash{\SetFigFont{12}{14.4}{\rmdefault}
{\mddefault}{\updefault}$\iota$}}}
\put(7189,-5561){\makebox(0,0)[lb]{\smash{\SetFigFont{12}{14.4}{\rmdefault}
{\mddefault}{\updefault}$\iota$}}}
\put(1201,-4100){\makebox(0,0)[lb]{\smash{\SetFigFont{12}{14.4}{\rmdefault}
{\mddefault}{\updefault}$\sigma$}}}
\put(9351,-4100){\makebox(0,0)[lb]{\smash{\SetFigFont{12}{14.4}{\rmdefault}
{\mddefault}{\updefault}$\sigma$}}}
\put(3401,-4100){\makebox(0,0)[lb]{\smash{\SetFigFont{12}{14.4}{\rmdefault}
{\mddefault}{\updefault}$\D'$}}}
\put(7101,-4100){\makebox(0,0)[lb]{\smash{\SetFigFont{12}{14.4}{\rmdefault}
{\mddefault}{\updefault}$\D'$}}}
\end{picture}
\end{center}

\vskip24pt

\section{Association} 
The concept of association or coupled calculus is one of the distinguishing 
features of local Colombeau algebras. It introduces a (linear) equivalence 
relation
on the algebra identifying those elements which are ``equal in the sense of
distributions'', thereby allowing to identify  ``distributional shadows'' of
certain elements of the algebra. 
This construction  amounts to determining the macroscopic
aspect of the regularization procedure encoded in them. 
Phrased more technically, a linear quotient of
$\ce_m$ resp.\ $\cg$ is formed containing $\D'$ as a subspace. Especially in
physical modelling this notion provides a useful tool for analyzing nonlinear
problems involving singularities (cf.\  e.g.\  \cite{CVW}, \cite{KS}, \cite{VW},
\cite{Vesi}). In what follows we extend the notion of association to
$\hat\G(M)$.
\bd
An element $[R]$ of $\hat\G(M)$ is called associated to $0$ $([R]\approx 0)$ if for
some (hence every) representative $R$ of $[R]$ we have: $\forall \om \in \Om_c^n(M)$
$\exists m>0$ with
\begin{equation} \label{ass}
\lim_{\eps\to 0} \int\limits_M R(\Phi(\eps,p),p)\om(p) = 0 \qquad \forall \Phi\in 
\tilde\mathcal{A}_m(M)
\end{equation}
Two elements $[R]$, $[S]$ of $\hat\G(M)$ are called associated $([R]\approx [S])$ if
$[R-S] \approx 0$. We say that $[R]\in \hat\G(M)$ admits $u\in \D'(M)$ as an associated
distribution if $[R]\approx \iota(u)$, i.e. if $\forall \om \in \Om_c^n(M)$
$\exists m>0$ with
\begin{equation} \label{assdist}
\lim_{\eps\to 0} \int\limits_M R(\Phi(\eps,p),p)\om(p) = \langle u, \om\rangle \qquad 
\forall \Phi\in \tilde\mathcal{A}_m(M)
\end{equation}
\et
Finally, by the same methods as in the local theory
we obtain consistency in the sense of association of
classical multiplication operations with multiplication in the algebra:
\bp \label{cons}
\begin{itemize}
  \item[(i)] If $f\in \cc^\infty(M)$ and $u\in \D'(M)$ then
  \begin{equation} \label{fuass}
    \iota(f)\iota(u) \approx \iota(fu)
  \end{equation} 
  \item[(ii)] If $f,g \in \cc(M)$ then
  \begin{equation} \label{fgass}
    \iota(f)\iota(g) \approx \iota(fg)
  \end{equation} 
\end{itemize}
\et

{\small
{\it Electronic Mail:}

\begin{tabular}{ll}
M.G.: & {\tt michael@mat.univie.ac.at}\\
M.K.: & {\tt Michael.Kunzinger@univie.ac.at}\\
R.S.: & {\tt Roland.Steinbauer@univie.ac.at}\\
J.V.: & {\tt jav@maths.soton.ac.uk}
\end{tabular}
}
\end{document}